\documentclass[a4paper,12pt,centertags,psamsfonts]{amsart}
\usepackage{amsmath,latexsym}  %
\usepackage{amssymb}
\usepackage{times}
\chardef\bslash=`\\ 

\newtheorem*{thm}{Theorem} 

\newtheorem*{prop}{Proposition}

\theoremstyle{definition}
\newtheorem{defn}{Definition}[]

\theoremstyle{remark}
\newtheorem*{rem}{Remark} 

\newcommand{\Na}{\mathbb{N}}

\newcommand{\Zi}{\mathbb{Z}}

\newcommand{\R}{\mathbb{R}}
\newcommand{\fR}{\mathfrak{R}}

\newcommand{\C}{\mathbb{C}}

\let\bsy\boldsymbol
\def\P{\bsy P}
\def\T{\bsy T}
\def\t{\bsy t}

\def\p{\bsy p}
\def\f{\bsy f}
\def\K{\bsy K}
\def\k{\bsy k}

\def\F{\bsy F}

\def\H{\mathcal{H}}
\newcommand{\A}{\mathcal{A}}
\newcommand{\M}{\mathcal{M}}

\newcommand{\Ss}{\mathcal{S}}

\newcommand{\wt}{\widetilde}

\DeclareMathOperator{\sgn}{sign}

\newcommand{\eval}[2][\right]{\relax
  \ifx#1\right\relax \left.\fi#2#1\rvert}

\let\abs=\envert

\let\norm=\enVert

\let\hnorm=\henVert
\newcommand{\kskob}[1]{\left[#1\right]}
\let\ksk=\kskob

\newcommand{\skoba}[1]{\left\langle#1\right\rangle}
\let\ska=\skoba

\let\hska=\hskoba

\let\krsk=\krskoba

\let\sets=\fskoba
\cleardoublepage 
\begin{document}
\renewcommand{\sectionmark}[1]{}
\title[unitary equivalence to
bi-Carleman operators]{Simultaneous unitary equivalence to
bi-Carleman operators with arbitrarily smooth kernels of Mercer type}
\author[I. M. Novitski\u i]{Igor M. Novitski\u i}
\address{Institute for Applied Mathematics, Russian Academy of Sciences,
92, Zaparina Street, Khabarovsk 680 000, Russia}
\email{novim@iam.khv.ru}
\thanks{Research supported in part by grant N 03-1-0-01-009 from
the Far-Eastern Branch of the Russian Academy of Sciences. This paper
was written in November 2003, when the author enjoyed the hospitality of the
Mathematical Institute of Friedrich-Schiller-University, Jena, Germany}
\keywords{Integral linear operator, bi-Carleman operator,
Hilbert-Schmidt operator, Carleman kernel, essential spectrum,
Lemari\'e-Meyer wavelet}
\subjclass[2000]{Primary 47B38, 47G10; Secondary 45P05}
\begin{abstract}
In this paper, we characterize the families of those bounded linear operators on
a separable Hilbert space which are simultaneously unitarily equivalent
to integral bi-Carleman operators on $L_2(\R)$ having
\textit{arbitrarily} smooth kernels of Mercer type.
The main result
is a qualitative sharpening of an earlier result of \cite{Nov:Lon}.
\end{abstract}
\maketitle

\section{Introduction. Main Result}
Throughout,  $\H$ will denote
a separable Hilbert space with the inner product
$\hska{\cdot,\cdot}$ and the norm $\hnorm{\cdot}$, $\fR(\H)$ the algebra of all bounded linear
operators on $\H$, and $\C$, and $\Na$, and $\Zi$, the complex
plane, the set of all positive integers, the set of all integers,
respectively. For an operator $A$ in $\fR(\H)$, $A^*$ will denote the
Hilbert space adjoint of $A$ in $\fR(\H)$.
Given an operator $T\in \fR(\H)$, define an operator set
\begin{equation*}\label{MT}
\M(T)=\krsk{T\fR(\H)\cup T^*\fR(\H)}\cap\krsk{\fR(\H)T^*\cup \fR(\H)T)},
\end{equation*}
where $S\fR(\H)$, $\fR(\H)S$ stand for the sets
$$
\sets{SA\mid A\in \fR(\H)},\quad\sets{AS\mid A\in \fR(\H)},$$
respectively.

Throughout,  $C(X,B)$, where $B$ is a Banach space (with norm
$\norm{\cdot}_B$), denote the Banach space (with the norm $\norm{f}_{C(X,B)}
=\sup\limits_{x\in X}\,\norm{f(x)}_B$) of
continuous $B$-valued functions defined on a
locally compact space $X$ and \textit{vanishing at infinity\/} (that is,
given any $f \in C(X,B)$ and $\varepsilon>0$, there exists a compact subset
$X(\varepsilon,f) \subset X$ such that $\norm{f(x)}_{B} < \varepsilon$
whenever $x \not\in  X(\varepsilon,f)$).

Let $\R$ be the real line $(-\infty,+\infty)$ with the Lebesgue measure,
and let $L_2=L_2(\R)$ be the Hilbert
space of (equivalence classes of) measurable complex-valued functions on
$\R$ equipped with the inner product
$$
\ska{f,g}=\int_\R f(s)\overline{g(s)}\,ds
$$
and the norm
$\norm{f}=\ska{f,f}^{\frac{1}2}$.

A linear operator
                 $T : L_2 \to  L_2$
is said to be \textit{integral\/}
if there exists a measurable function
$\T$ on the Cartesian product $\R^2=\R\times\R$, a \textit{kernel\/},
such that, for every $f\in L_2$,
$$
               (Tf)(s)=\int_{\R} \T(s,t)f(t)\,dt
$$
for almost every $s$ in $\R$. A kernel $\T$ on
$\R^2$ is said to be \textit{Carleman\/} if $\T(s,\cdot) \in L_2$
for almost every fixed $s$ in $\R$.
An integral operator with a kernel $\T$ is
called \textit{Carleman\/} if $\T$ is a Carleman kernel, and it is called
\textit{bi-Carleman\/} if both $\T$ and $\T^*$ ($\T^*(s,t)=\overline{\T(t,s)}$)
are Carleman kernels.
Every Carleman kernel, $\T$, induces a \textit{Carleman
function\/} $\t$ from $\R$ to $L_2$ by
$\t(s)=\overline{\T(s,\cdot)}$
for all $s$ in ${\R}$ for which $\T(s,\cdot)\in L_2$.

We shall also recall a characterization of bi-Carleman
representable operators. Its version for self-adjoint operators was
first obtained by von Neumann \cite{Neu} and was later extended by Korotkov
to the general case (see \cite[p.\ 100]{Kor:book1},
\cite[p.\ 103]{Halmos:Sun}). The assertion says
that a necessary and sufficient condition that an operator
$S\in \fR(\H)$ be
unitarily equivalent to a bi-Carleman operator is that there exist an
orthonormal sequence $\left\{e_n\right\}$ such that
\begin{equation}\label{kh}
\hnorm{Se_n}\rightarrow 0,\quad
\hnorm{S^*e_n}\rightarrow 0\quad
\text{as $n\rightarrow\infty$}
\end{equation}
(or, equivalently, that $0$ belong to the essential spectrum of $SS^*+S^*S$).

\begin{defn}
Given any non-negative integer $m$, we say that a
function $\K$ on $\R^2$ is a $K^m$-{\it kernel\/}
(see \cite{Nov:Lon}, \cite{nov:91}) if
\begin{enumerate}
\renewcommand{\labelenumi}{(\roman{enumi})}
\item the function $\K$ and all its partial derivatives
on $\R^2$ up to order $m$ are in $C(\R^2,\C)$,
\item the Carleman function $\k$,
$\k(s)=\overline{\K(s,\cdot)}$,
and all its (strong) derivatives on $\R$ up to order $m$
are in $C(\R,L_2)$,
\item the conjugate transpose function $\K^*$,
$\K^*(s,t)=\overline{\K(t,s)}$,
satisfies Condition (ii), that is,
the Carleman function $\k^*$, $\k^*(s)
=\overline{\K^*(s,\cdot)}$, and all its (strong) derivatives on $\R$ up
to order $m$ are in $C(\R,L_2)$.
\end{enumerate}
In addition, we say  that a function $\K$ is a
\textit{$K^\infty$-kernel} (see \cite{Je}, \cite{SP}) if it is a
$K^m$-kernel for each non-negative integer $m$.
\end{defn}

\begin{defn}
Let $\K$ be a $K^m$($K^\infty$)-kernel and let $T$ be the integral operator
it induces. We say that the $K^m$($K^\infty$)-kernel $\K$ is of
{\it Mercer type\/} if every operator $A\in\M(T)$ is an
integral operator having $K^m$($K^\infty$)-kernel.
\end{defn}

The concept of Mercer type $K^m$-kernels for finite $m$ was first introduced
in our paper \cite{Nov:Lon} where there is a motivation of the reason why this subclass
of $K^m$-kernels deserves the qualification ``of Mercer type''.

Given any non-negative integer $m$, the following result both
gives a characterization of all bounded operators whose unitary orbits
contain a bi-Carleman operator having $K^m$-kernel of Mercer type
and describes families of those operators that can be
simultaneously unitarily represented as bi-Carleman operators having
$K^m$-kernels of Mercer type
(cf. \eqref{kh}).

\begin{prop}[\cite{Nov:Lon}]\label{msmooth}
If for an operator family $\left\{S_\alpha\mid\alpha \in \A\right\}\subset \fR(\H)$
there exists an orthonormal sequence
$\left\{e_n\right\}$
such that
$$
\lim\limits_{n\to\infty}\sup\limits_{\alpha\in\A}\hnorm{S_\alpha^*e_n}=0,
\quad
\lim_{n\to\infty}\sup_{\alpha\in\A}\hnorm{S_\alpha e_n}=0,
$$
then there exists a unitary operator $U_m:\H\to L_2$
such that all the operators $U_mS_\alpha U_m^{-1}$ $(\alpha\in\A)$
and their linear combinations are bi-Carleman operators having $K^m$-kernels of Mercer type.
\end{prop}
The construction of the unitary operator $U_m$ given in the proof
of Proposition  depends on the preassigned order $m<\infty$ of smoothness
(see \cite{Nov:Lon}). The purpose of the present paper is to show that
Proposition is true with $K^\infty$-kernels in the conclusion, that is, to
prove the following qualitative sharpening of Proposition.

\begin{thm}\label{infsmooth}
If for an operator family $\sets{S_\alpha\mid\alpha \in \A}\subset \fR(\H)$
there exists an orthonormal sequence
$\sets{e_n}$
such that
\begin{equation}\label{1.2}
\lim\limits_{n\to\infty}\sup\limits_{\alpha\in\A}\hnorm{S_\alpha^*e_n}=0,
\quad
\lim_{n\to\infty}\sup_{\alpha\in\A}\hnorm{S_\alpha e_n}=0,
\end{equation}
then there exists a unitary operator
$U_\infty:\H\to L_2$ such that all the operators
$U_\infty S_\alpha U_\infty^{-1}$ $(\alpha\in\A)$
and their linear combinations
are bi-Carleman operators
having $K^\infty$-kernels of Mercer type.
\end{thm}

\section{Proof of Theorem}
The proof is broken up into three steps. The first step is to find
suitable orthonormal bases
$\sets{u_n}$ in $L_2$ and $\sets{f_n}$ in
$\H$ on which the construction of $U_\infty$ will be based.
The next step is to define a certain unitary operator that
sends the basis $\sets{f_n}$ onto the basis $\sets{u_n}$.
This operator is suggested as $U_\infty$ in the theorem, and
the rest of the proof is a straightforward verification that it is
indeed as desired. Thus, the proof yields more than just existence of
the unitary equivalence; it yields an explicit construction of the
unitary operator. From the point of view of the applications to operator
equations, the explicit computability of $U_\infty$ is an important side issue.
\subsection*{Step 1.}
For the proof, it will be convenient to have the following notation: if
an equivalence class $f\in L_2$ contains a function belonging to
$C(\R,\C)$, then we shall use
$\ksk{f}$ to denote that function.

Let $\left\{S_\alpha\mid \alpha\in\A\right\}\subset\fR(\H)$ be a family
satisfying \eqref{1.2} with the orthonormal
sequence $\sets{e_n}_{n=1}^\infty$.
Take orthonormal bases $\sets{f_n}$ for $\H$ and $\sets{u_n}$ for $L_2$
 which satisfy the conditions:
\begin{enumerate}
\renewcommand{\labelenumi}{(\alph{enumi})}
\item the terms of the sequence
$\left\{\ksk{u_n} ^{(i)}\right\}$ of derivatives are in $C(\R,\C)$, for each $i$
(here and throughout, the letter $i$ is reserved for integers in $[0,+\infty)$),
\item $\{u_n\}=\{g_k\}_{k=1}^\infty\cup\{h_k\}_{k=1}^\infty$, where
$\{g_k\}_{k=1}^\infty\cap\{h_k\}_{k=1}^\infty=\varnothing$,
and,  for each $i$,
\begin{equation}\label{hki}
\sum_k H_{k,i}<\infty\quad
\text{with $H_{k,i}=\norm{\ksk{h_k}^{(i)}}_{C(\R,\C)}$}\quad (k\in\Na)
\end{equation}
(the sum notation $\sum\limits_k$ will always be used instead of the more detailed
symbol $\sum\limits_{k=1}^\infty$),
\item $\sets{f_n}=\sets{x_k}_{k=1}^\infty\cup\sets{y_k}_{k=1}^\infty$
where $\sets{x_k}_{k=1}^\infty\cap\sets{y_k}_{k=1}^\infty=\varnothing$,
$\sets{x_k}_{k=1}^\infty\subset\sets{e_n}_{n=1}^\infty$,
and, for each $i$,
\begin{equation}\label{zndn}
\sum_k d_k\krsk{G_{k,i}+1}<\infty
\end{equation}
with $d_k=2\left(\sup\limits_\alpha\hnorm{S_\alpha x_k}^\frac{1}4+
\sup\limits_\alpha\hnorm{S_\alpha^*x_k}^\frac{1}4\right)\le1$,
and $G_{k,i}=\norm{\ksk{g_k}^{(i)}}_{C(\R,\C)}$ ($k\in\Na$).
\end{enumerate}
The proof uses the bases just described to construct the desired unitary
operator $U_\infty$.
\begin{rem}
Let $\sets{u_n}$ be an orthonormal basis for $L_2$ such that, for
each $i$,
\begin{gather}\label{1}
\kskob{u_n}^{(i)}\in C(\R,\C)\quad(n\in\Na),\\
\label{2}\norm{\kskob{u_n}^{(i)}}_{C(\R,\C)}\le D_nA_i\quad(n\in\Na),\\
\label{3}
\sum_kD_{n_k}<\infty,
\end{gather}
where $\{D_n\}_{n=1}^\infty$, $\{A_i\}_{i=0}^\infty$ are sequences
of positive numbers, and $\sets{n_k}_{k=1}^\infty$ is a subsequence
of $\Na$ such that
$\Na\setminus\sets{n_k}_{k=1}^\infty$ is a countable set.
Since $d(e_n)\to0$ as $n\to\infty$, it follows
that there exists a subset
$\sets{x_k}_{k=1}^\infty\subset\sets{e_n}_{n=1}^\infty$
for which Condition \eqref{zndn} holds with
$\sets{g_k}_{k=1}^\infty=\sets{u_n}\setminus\sets{u_{n_k}}_{k=1}^\infty$.
Moreover, the properties \eqref{2} and \eqref{3} imply Condition \eqref{hki} for
$h_k=u_{n_k}$ ($k\in\Na$).
Complete the set $\sets{x_k}_{k=1}^\infty$ to an orthonormal basis,
and let $y_k$ ($k\in\Na$) denote the new elements of that basis.
Then the bases $\sets{f_n}=\sets{x_k}_{k=1}^\infty\cup\sets{y_k}_{k=1}^\infty$
and $\{u_n\}$ satisfy Conditions (a)-(c).

A good example of the basis satisfying \eqref{1}-\eqref{3}
is a basis generated by the Lemari\'e-Meyer wavelet
$$
u(s)=\dfrac1{2\pi}\int_{\R}e^{i\xi(\frac12+s)}
\sgn\xi b(|\xi|)\,d\xi\quad (s\in\R),
$$
with the bell function $b$ belonging to
$C^\infty(\R)$ (for construction of the Lemari\'e-Meyer wavelets
we refer to \cite{LeMe}, \cite[\S~4]{Ausch}, \cite[Example D, p.~62]{Her}).
In this case, $u$ belongs to the Schwartz class $\Ss(\R)$, and hence
all the derivatives $\kskob{u}^{(i)}$
are in $C(\R,\C)$.
The corresponding orthonormal basis for $L_2$ is given by
$$
u_{jk}(s)=2^{\frac j2}u(2^js-k)\quad  (j,\,k\in\Zi).
$$
 Rearrange, in a completely arbitrary manner, the orthonormal set
$\{u_{jk}\}_{j,\,k\in\Zi}$ into a simple sequence,
so that it becomes $\{ u_n\}_{n\in\Na}$. Since, in view of this rearrangement,
to each $n\in\Na$ there corresponds a unique pair of integers
$j_n$, $k_n$, and
conversely, we can write, for each $i$,
$$
\norm{\kskob{u_n}^{(i)}}_{C(\R,\C)}=\norm{\kskob{u_{j_nk_n}}^{(i)}}_{C(\R,\C)}
\le D_nA_i,
$$
where
$$
D_n=\begin{cases}
2^{j_n^2}&\text{if $j_n>0$,}\\
\left(\dfrac1{\sqrt{2}}\right)^{\abs{j_n}}&\text{if $j_n\le0$,}
\end{cases}
\qquad A_i=2^{\left(i+\frac12\right)^2}\norm{\kskob{u}^{(i)}}_{C(\R,\C)}.
$$
Whence it follows that if $\{n_k\}_{k=1}^\infty\subset\Na$ is a subsequence
such that $j_{n_k}\to -\infty$ as $k\to\infty$, then
$$\sum_kD_{n_k}<\infty.$$
Thus,  the basis $\{u_n\}$ satisfies Conditions \eqref{1}-\eqref{3}.
\end{rem}

\subsection*{Step 2.} In this step our intention is to construct a candidate
for the desired unitary operator $U_\infty$ in the theorem.
Define such a unitary operator $U_\infty:\H\to L_2$ on the basis vectors by
setting
\begin{equation}\label{2.12}
U_\infty x_k=g_k,\quad U_\infty y_k=h_k
\quad\text{for all $k\in \Na$,}
\end{equation}
in the harmless assumption that $U_\infty f_n=u_n$ for all $n\in\Na$.

\subsection*{Step 3.}
    The verification that $U_\infty$ in \eqref{2.12} has the desired properties
is straightforward. Fix an arbitrary $\alpha\in\A$ and put
$T=U_\infty S_\alpha U_\infty ^{-1}$. Once this is done,
the index $\alpha$ may be omitted for $S_\alpha$.

Let $E$ be the orthogonal projection onto
the closed linear span of the vectors $x_k$ ($k\in\Na$).
Split the operator $S$ as follows:
\begin{equation}\label{2.7}
S=(1-E)S+ES,\quad S^*=(1-E)S^*+ES^*.
\end{equation}
The operators $J=S E$ and
$\wt J=S^* E$ are nuclear operators and,
therefore, are Hilbert--Schmidt operators; these properties are almost
immediate consequences of \eqref{zndn}.
\par
Write the Schmidt decompositions
$$
J=\sum_n s_{n}\hska{\cdot,p_{n}} q_{n},
\quad\wt J=\sum_n \wt s_{n}\hska{\cdot,\wt p_{n}
}\wt q_{n},
$$
where the
$s_{n}$
are the singular values of $J$
(eigenvalues of $\left(J J^*\right)^{\frac1{2}}$),
$\sets{p_n}$, $\sets{q_n}$  are
orthonormal sets (the $p_{n}$ are eigenvectors for
$J^* J$ and $q_{n}$ are eigenvectors for
$J J^*$).
The explanation of the notation for $\wt J$ is similar.
\par
Now introduce auxiliary operators
$B$, $\wt B$
by
\begin{equation}\label{aux}
B=\sum_n s_{n}^{\frac1{4}}
\hska{\cdot,p_{n}} q_{n},
\quad
\wt B=\sum_n \wt s_{n}^{\frac1{4}}\hska{\cdot,
\wt p_{n}}\wt q_{n}.
\end{equation}
The Schwarz inequality yields
\begin{equation}\label{2.8}
\begin{gathered}
\|B^* x_k\|_\H+
\hnorm{B x_k}+
\hnorm{\wt B^* x_k}+
\hnorm{\wt B x_k}
\\
=
\hnorm{\left(J J^*\right)^{\frac1{8}}x_k}+
\hnorm{\left(J^* J\right)^{\frac1{8}} x_k}\\+
\hnorm{\left(\wt J\wt J^*\right)^{\frac1{8}}x_k}+
\hnorm{\left(\wt J^*\wt J\right)^{\frac1{8}} x_k}
\\
\leq
\hnorm{J^* x_k}^{\frac1{4}}+
\hnorm{J x_k}^{\frac1{4}}+
\hnorm{\wt J^* x_k}^{\frac1{4}}+
\hnorm{\wt J x_k}^{\frac1{4}}
\leq d_k.
\end{gathered}
\end{equation}
It follows that all the operators $B$, $\wt B$ are nuclear
operators (see \eqref{zndn}) and hence
\begin{equation}\label{2.9}
\sum_n s_{n}^{\frac1{2}}<\infty,\quad
\sum_n \wt s_{n}^{\frac1{2}}<\infty.
\end{equation}
Define
$Q=(1-E)S^*$, $\wt Q=(1-E)S$.
Then Condition (c) provides the representations
\begin{equation}
\begin{gathered}\label{2.10}
Q f=\sum_k \hska{ Q f,y_k} y_k=
          \sum_k \hska{ f,S y_k} y_k,\\
\wt Q f=\sum_k \hska{\wt Q f,y_k} y_k=
          \sum_k \hska{ f,S^* y_k} y_k,
\end{gathered}
\end{equation}
for all $f$ in $\H$.
\par

Using the decompositions \eqref{2.7}, which now look like
$S=\wt Q+\wt J^*$,
$S^*=Q+J^*$,
we shall prove presently that $T$ is an integral
operator having $K^\infty$-kernel of Mercer type.
\par
From \eqref{2.10} and \eqref{2.12}, it follows that, for each $f\in L_2$,
\begin{equation}\label{P}
\begin{gathered}
Pf=U_\infty QU_\infty^{-1}f=\sum_k \hska{f,Th_k} h_k,\\
\wt Pf=U_\infty \wt QU_\infty ^{-1}f=
\sum_k\hska{f,T^*h_k} h_k.
\end{gathered}
\end{equation}
Represent the equivalence classes
$Th_k$, $T^*h_k$ ($k\in\Na$) by the Fourier
expansions
$$
Th_k=\sum_n \hska{ y_k,S^*f_n}u_n,
\quad
T^*h_k=\sum_n \hska{ y_k,Sf_n} u_n,
$$
where the series converge in the $L_2$ sense. But more than that can be
said about convergence, namely that, for each fixed $i$, the series
\begin{equation}\label{2.15}
\sum_n \hska{y_k,S^*f_n} \ksk{u_n}^{(i)}(s),
\quad
\sum_n \hska{y_k,Sf_n} \ksk{u_n}^{(i)}(s)\quad(k\in\Na)
\end{equation}
converge in the norm of $C(\R,\C)$.
Indeed, all the series are everywhere pointwise dominated by one series
$$
\sum_n\left(\hnorm{S^*f_n}+\hnorm{Sf_n}\right)\abs{ \ksk{u_n}^{(i)}(s)},
$$
which is uniformly convergent on $\R$ for the following reason: its
subseries
$$
\begin{gathered}
\sum_k\left(\hnorm{Sx_k}+\hnorm{S^*x_k}\right)\abs{\ksk{g_k}^{(i)}(s)},\\
\sum_k\left(\hnorm{Sy_k}+\hnorm{S^*y_k}\right)\abs{\ksk{h_k}^{(i)}(s)}
\end{gathered}
$$
are uniformly convergent on $\R$ because they in turn are dominated by the
convergent series
\begin{equation}\label{2.16}
\sum_k d_kG_{k,i},
\quad
\sum_k 2\|S\|H_{k,i},
\end{equation}
respectively (see \eqref{zndn}, \eqref{hki}).
\par
It is now evident that the pointwise sums in \eqref{2.15} define
functions that belong to $C(\R,\C)$.
Moreover, the above arguments prove that, for each fixed $i$, the
derivative sequences
$\sets{\ksk{ Th_k}^{(i)}}$,
$\sets{\ksk{ T^*h_k}^{(i)}}$
are uniformly bounded in $C(\R,\C)$ in the sense
that there exists a positive constant $C_i$ such that
$$
\norm{\ksk{ Th_k}^{(i)}}_{C(\R,\C)}<C_i,
\quad
\norm{\ksk{ T^*h_k}^{(i)}}_{C(\R,\C)}<C_i,
$$
for all $k$. Hence, by \eqref{hki},
it is possible to infer that, for all non-negative integers
$i$, $j$, both
$$
\sum_k \ksk{h_k}^{(i)}(s)
\overline{\ksk{ Th_k}^{(j)}(t)}\quad
\text{and}\quad
\sum_k \ksk{h_k}^{(i)}(s)
\overline{\ksk{ T^*h_k}^{(j)}(t)}
$$
converge in the norm of $C(\R^2,\C)$.
This makes it obvious that both
\begin{equation}\label{2.14}
\begin{gathered}
\P(s,t)=\sum_k \ksk{ h_k}(s)\overline
{\ksk{ Th_k}(t)}\\
\text{and}\\
\wt \P(s,t)=\sum_k \ksk{ h_k}(s)\overline{\ksk{ T^*h_k}(t)},
\end{gathered}
\end{equation}
satisfy Condition (i) for each $m$.
\par
Now we prove that the (Carleman) functions
\begin{equation}\label{pp}
\begin{gathered}
\p(s)=\overline{\P(s,\cdot)}=\sum_k\overline{\ksk{h_k}(s)}Th_k,\\
\wt \p(s)=\overline{\wt \P(s,\cdot)}=
\sum_k\overline{\ksk{h_k}(s)}T^*h_k
\end{gathered}
\end{equation}
satisfy Condition (ii) for all $m$.
Indeed, the series displayed converge absolutely in the
$C(\R,L_2)$ sense, because those two series whose terms are
$\abs{ \ksk{h_k}(s)}\norm{Th_k}$
and
$\abs{ \ksk{h_k}(s)}\norm{T^*h_k}$
respectively are dominated by
the second series in \eqref{2.16} for $i=0$.
For the remaining $i$, a similar
reasoning implies the same conclusion for the series
$$
\sum_k\overline{\ksk{h_k}^{(i)}(s)}Th_k,
\quad
\sum_k\overline{\ksk{h_k}^{(i)}(s)}T^*h_k.
$$
The asserted property of both $\p$ and $\wt \p$
to satisfy (ii) for each $m$ then follows from the termwise differentiation
theorem.
Now observe that, by \eqref{hki} and \eqref{pp}, the series in \eqref{P}
(viewed, of course, as ones with terms belonging to $C(\R,\C)$) converge
(absolutely) in $C(\R,\C)$-norm  to the functions
$$
\begin{gathered}
\kskob{Pf}(s)\equiv\skoba{f,\p(s)}\equiv\int_{\R}\P(s,t)f(t)\,dt,\\
\kskob{\wt Pf}(s)\equiv\skoba{f,\wt\p(s)}\equiv\int_{\R}
\wt\P(s,t)f(t)\,dt,
\end{gathered}
$$
respectively.
Thus, both $P$ and $\wt P$ are Carleman operators
with $\P$ and $\wt \P$ their kernels, respectively, satisfying Conditions
(i), (ii) for each $m$.

Now consider the (integral) Hilbert--Schmidt operators
$F=U_\infty J^*U_\infty ^{-1}$ and $\wt F=U_\infty \wt J^*U_\infty ^{-1}$.
Prove that both $F$ and $\wt F$ have kernels satisfying
(i) for each $m$.
Starting from the Schmidt decompositions for $F$ and $\wt F$, define their
kernels by
\begin{equation}\label{2.17}
\begin{gathered}
\F(s,t)
=\sum_n s_n^{\frac1{2}}\ksk{ U_\infty B^*q_n} (s)
\overline{\ksk{ U_\infty Bp_n} (t)},\\
\wt \F(s,t)
=\sum_n\wt s_n^{\frac1{2}}\ksk{ U_\infty \wt B^*\wt q_n
}(s)\overline{\ksk{ U_\infty \wt B\wt p_n}(t)},
\end{gathered}
\end{equation}
for all $s$, $t$ in $\R$, in the tacit assumption that the square brackets
are everywhere permissible (for the auxiliary operators $B$, $\wt B$
see \eqref{aux}). In view of \eqref{2.9} the desired conclusion that the kernels
so defined satisfy (i) for each $m$ can be inferred as soon as it is known that for each
fixed $i$ the terms of the sequences
$$
\sets{\ksk{ U_\infty Bp_k}^{(i)}},
\quad
\sets{\ksk{ U_\infty B^*q_k}^{(i)}},
\quad
\sets{\ksk{ U_\infty \wt B\wt p_k}^{(i)}},
\quad
\sets{\ksk{ U_\infty \wt B^*\wt q_k}^{(i)}}
$$
make sense, are in $C(\R,\C)$, and are uniformly bounded in $C(\R,\C)$.

To see the validity of the properties indicated, observe that all the series
\begin{equation*}
\begin{gathered}
\sum_n\hska{p_k,B^*f_n} \ksk{u_n}^{(i)}(s),\quad
\sum_n\hska{q_k,Bf_n} \ksk{u_n}^{(i)}(s),\quad
\\\sum_n\hska{\wt p_k,\wt B^*f_n} \ksk{u_n}^{(i)}(s),\quad
\sum_n\hska{\wt q_k,\wt Bf_n} \ksk{u_n}^{(i)}(s)
\quad(k\in\Na)
\end{gathered}
\end{equation*}
(which in the case where $i=0$ are just the Fourier expansions for
$U_\infty Bp_k$, $U_\infty B^*q_k$, $U_\infty \wt B\wt p_k$, $U_\infty \wt B^*\wt q_k$)
are dominated by one series
$$
\sum_nc(f_n)\abs{ \ksk{u_n}^{(i)}(s)},
$$
where
$c(g)=\hnorm{B^*g}+\hnorm{Bg}+\hnorm{\wt B^*g}+\hnorm{\wt Bg}$
whenever
$g\in \H$. The last series is uniformly convergent, because
it consists of the two dominatedly and uniformly convergent subseries
$$
\sum_nc(x_k)\abs{ \ksk{g_k}^{(i)}(s)},
\quad
\sum_nc(y_k)\abs{ \ksk{h_k}^{(i)}(s)};
$$
the corresponding dominant series are
$$
\sum_kd_kG_{k,i},
\quad
\sum_k2\left(\|B\|+\|\wt B\|\right)H_{k,i}
$$
(see \eqref{2.8}, \eqref{zndn}, \eqref{hki}).
\par
In view of  \eqref{2.9} and the uniform boundedness in $C(\R,\C)$ of the
sequences
$\sets{\ksk{ U_\infty B^*q_n}^{(i)}}$,
$\sets{\ksk{ U_\infty \wt B^*\wt q_n}^{(i)}}$
for each fixed $i$, the series
$$
\sum_n s_n^{\frac1{2}}\overline{\ksk{ U_\infty B^*q_n}^{(i)}(s)}U_\infty Bp_n,
\quad
\sum_n\wt s_n^{\frac1{2}}\overline{\ksk{ U_\infty \wt B^*\wt q_n
}^{(i)}(s)}U_\infty \wt B\wt p_n
$$
are absolutely convergent in the $C(\R,L_2)$ sense, and hence their sums
belong to $C(\R,L_2)$. Observe by \eqref{2.17} that two of them, namely those
for $i=0$, represent the Carleman  functions
$\f(s)=\overline{\F(s,\cdot)}$,
$\wt \f(s)=\overline{\wt \F(s,\cdot)}$.
Thus, both Carleman functions $\f$ and $\wt \f$ satisfy
Condition (ii) for every $m$.
\par

     In accordance with  \eqref{2.7}, the operator $T$, which is the transform
by $U_\infty $ of $S$, has the decompositions
$T=\wt P+\wt F$,
$T^*=P+F$
where all the terms are the Carleman operators already described.
So both $T$ and $T^*$ are Carleman operators, and their kernels
$\K$ and $\wt \K$, which are defined by
\begin{equation}\label{2.18}
\K(s,t)=\wt \P(s,t)+\wt \F(s,t),
\quad
\wt \K(s,t)=\P(s,t)+\F(s,t),
\end{equation}
for all $s$, $t\in\R$, inherit the two properties (i), (ii) from their terms,
for each $m$.
Since (cf. \cite[p.\ 37]{Halmos:Sun}) it is possible to write
$\K(s,t)=\overline{\wt \K(t,s)}$ and
$\K(\cdot,t)=\overline{\wt \K(t,\cdot)}$ for all $s$, $t\in \R$, the kernel
$\K$ satisfies Conditions (i), (ii), (iii) for each $m$, so that it is a
$\K^\infty$-kernel.
\par
As the preceding proof shows, the major condition that an operator
$A\in \fR(\H)$
must satisfy in order that $U_\infty AU_\infty ^{-1}$
be an integral operator having $K^\infty$-kernel is that, for each
$k$, $$
2\left(\hnorm{Ax_k}^{\frac1{4}}+\hnorm{A^*x_k}^{\frac1{4}}\right)\leq
d_k. $$ If $A\in \M(S)$ then there exist operators
$V$, $W\in \fR(\H)$ such that at least one of the relations
$A=SV=WS$, $A=S^*V=WS^*$, $A=SV=WS^*$, $A=VS=S^*W$ holds. In any event,
whatever its decomposition may be, the operator $A$ satisfies the
inequalities $$
\left(\hnorm{Ax_k}^{\frac1{4}}+\hnorm{A^*x_k}^{\frac1{4}}\right)\leq
2c\left(\hnorm{Ax_k}^{\frac1{4}}+\hnorm{A^*x_k}^{\frac1{4}}\right)\leq
cd_k, $$ where $c^4=\max\sets{\|V\|,\|W\|}$. This implies, by the
above remark, that $U_\infty $ automatically carries every
$A\in\M(S)$ onto an integral operator $U_\infty AU_\infty ^{-1}$
having $K^\infty$-kernel so that $\K$ in  \eqref{2.18} is a
$K^\infty$-kernel of Mercer type.
\par
The fact that those $K^\infty$-kernels which induce finite linear
combinations of $U_\infty S_\alpha U_\infty ^{-1}$ are of Mercer type remains to be proved.
The result can be inferred from the result for $S$, which has just been obtained.
Indeed, consider any finite linear combination $G=\sum z_\alpha S_\alpha$
with $\sum \abs{ z_\alpha}\leq 1$.
It is seen easily that, for each $n$,
$$
\hnorm{\sum z_\alpha S_\alpha e_n}\leq\sup_\alpha\hnorm{S_\alpha e_n},
\quad
\hnorm{\sum\overline{z}_\alpha S^*_\alpha e_n}
\leq\sup_\alpha\hnorm{S^*_\alpha e_n}.
$$
There is, therefore, no barrier to assuming that $G$ was, from the start,
in $\sets{S_\alpha\mid\alpha\in\A}$ and even equal to $S$.
The proof of Theorem  is complete.

\section*{Acknowledgments} The author thanks the Mathematical Institute
of the  University of Jena for its hospitality, and specially
W.~Sickel and H.-J.~Schmei\ss er for useful remarks and fruitful
discussion on applying wavelets in integral representation theory.

\bibliographystyle{amsplain}

\end{document}